\documentclass[12pt]{amsart}

\usepackage{a4wide}

\newtheorem{theorem}{Theorem}[section]

\theoremstyle{definition}

\newcommand{\C}{\ensuremath{\mathbb{C}}}
\newcommand{\R}{\ensuremath{\mathbb{R}}}
\newcommand{\g}[1]{\ensuremath{\mathfrak{#1}}}

\DeclareMathOperator{\tr}{tr}

\begin{document}
\title[Inhomogeneous isoparametric hypersurfaces]{Inhomogeneous isoparametric
hypersurfaces\\ in complex hyperbolic spaces}

\author[J.C. D\'\i{}az-Ramos]{Jos\'{e} Carlos D\'\i{}az-Ramos}
\address{Department of Geometry and Topology,
University of Santiago de Compostela, Spain.}
\email{josecarlos.diaz@usc.es}

\author[M. Dom\'{\i}nguez-V\'{a}zquez]{Miguel Dom\'{\i}nguez-V\'{a}zquez}
\address{Department of Geometry and Topology,
University of Santiago de Compostela, Spain.}
\email{miguel.dominguez@usc.es}

\thanks{The first author has been supported by a Marie-Curie
European Reintegration Grant (PERG04-GA-2008-239162). The second
author has been supported by the FPU programme of the Spanish
Government. Both authors have been supported by projects
MTM2006-01432 and INCITE09207151PR (Spain).}



\begin{abstract}
We construct examples of inhomogeneous isoparametric real
hypersurfaces in complex hyperbolic spaces.
\end{abstract}


\maketitle

\section{Introduction}

An isoparametric hypersurface of a Riemannian manifold is a
hypersurface that is a level set of an isoparametric function. Cartan
proved that a hypersurface is isoparametric if and only if the
hypersurface itself and its sufficiently close parallel hypersurfaces
have constant mean curvature. If the ambient manifold has constant
sectional curvature, then Cartan also proved that a hypersurface is
isoparametric if and only if it has constant principal curvatures.
The study of isoparametric hypersurfaces is an active topic of
research in Differential Geometry. See~\cite{T00} for a survey.

In a more general ambient space, an isoparametric hypersurface might
have nonconstant principal curvatures, and therefore it might be
inhomogeneous. See~\cite{W82} for isoparametric hypersurfaces with
nonconstant principal curvatures in complex projective spaces. These
examples are constructed from the inhomogeneous examples in spheres
given by Ferus, Karcher and M\"{u}nzner in~\cite{FKM81} via the Hopf map.
To our knowledge, the only examples of isoparametric hypersurfaces in
complex hyperbolic spaces known so far were homogeneous. Homogeneous
hypersurfaces in complex hyperbolic spaces were classified
in~\cite{BT07} and their geometry was studied in~\cite{BD09}.

The aim of this paper is to construct examples of isoparametric
hypersurfaces in complex hyperbolic spaces that are in general not
homogeneous. These are not related to the Ferus, Karcher and M\"{u}nzner
hypersurfaces in~\cite{FKM81}. Our examples arise as tubes around
certain homogeneous submanifolds that are in a way a modification of
the homogeneous submanifolds introduced by Berndt and Br\"{u}ck
in~\cite{BB01}. Indeed, the tubes around Berndt-Br\"{u}ck submanifolds
are a particular case of our examples, and are the only ones being
homogeneous hypersurfaces. Xiao claims in~\cite{X99} that
isoparametric real hypersurfaces in $\C H^n$ are homogeneous,
although there is no proof of this fact to our knowledge.
Furthermore, our examples show that inhomogeneous isoparametric
examples do exist.

The main result of our paper, the construction of these isoparametric
hypersurfaces, is summarized in Theorem~\ref{thMain}. In general, our
examples are inhomogeneous. Pointwise, the principal curvatures of
these hypersurfaces are the same as those of the Berndt and Br\"{u}ck
examples. Then, we also show that the inhomogeneous hypersurfaces in
our examples correspond precisely to those that have nonconstant
principal curvatures.

\section{Preliminaries}\label{secPrelimnaries}

In this paper we follow the notation of~\cite{BD09}.

Let $\C H^n$ be the complex hyperbolic space of constant holomorphic
sectional curvature $-1$. We write $\C H^n$ as $G/K$ where
$G=SU(1,n)$ and $K=S(U(1)U(n))$. As usual, denote by $\g{g}$ and
$\g{k}$ the Lie algebras of $G$ and $K$, respectively. Let
$\g{g}=\g{k}\oplus\g{p}$ be the Cartan decomposition of $\g{g}$ with
respect to a point $o\in\C H^n$ and choose a maximal abelian subspace
$\g{a}$ of $\g{p}$, which is therefore $1$-dimensional. Let
$\g{g}=\g{g}_{-2\alpha}\oplus\g{g}_{-\alpha}\oplus\g{g}_{0}
\oplus\g{g}_{\alpha}\oplus\g{g}_{2\alpha}$ be the root space
decomposition of $\g{g}$ with respect to $o$ and $\g{a}$. We
introduce an ordering in the set of roots so that $\alpha$ is a
positive root. Then $o$, $\g{a}$, and this ordering determine a point
at infinity $x$ in the ideal boundary $\C H^n(\infty)$ of $\C H^n$.
Let $\g{n}=\g{g}_\alpha\oplus\g{g}_{2\alpha}$. Then
$\g{g}=\g{k}\oplus\g{a}\oplus\g{n}$ is the Iwasawa decomposition of
the Lie algebra $\g{g}$ with respect to the point $o\in\C H^n$ and
the point at infinity $x\in\C H^n(\infty)$. Let $A$, $N$ and $AN$ be
the connected and simply connected subgroups of $G$ whose Lie
algebras are $\g{a}$, $\g{n}$, and $\g{a}\oplus\g{n}$, respectively.
Thus $G$ is diffeomorphic to $K\times A \times N$, and $AN$ is
diffeomorphic to $\C H^n$. The metric (resp., the complex structure
$J$) on $\C H^n$ induces a left invariant metric (resp.\ a complex
structure $J$) on the solvable Lie group $AN$, so that $\C H^n$ and
$AN$ become isometric as K\"{a}hler manifolds. We also have the
isomorphism $T_o\C H^n\cong\g{a}\oplus\g{n}$.

From now on, let $B$ be the unit left-invariant vector field of
$\g{a}$ determined by $x$; that is, the geodesic through $o$ whose
initial speed is $B$ converges to the point at infinity $x$. Set
$Z=JB\in\g{g}_{2\alpha}$. We obviously have $\g{a}=\R B$ and
$\g{g}_{2\alpha}=\R Z$. Moreover, $\g{g}_\alpha$ is $J$-invariant, so
it is isomorphic to $\C^{n-1}$. The Lie algebra structure on
$\g{a}\oplus\g{n}$ is given by the relations
\[
\left[B,Z\right]=Z,\quad 2\left[B,U\right]=U,\quad
\left[U,V\right]=\langle JU, V\rangle Z,\quad \left[Z,U\right]=0,
\]
where $U$, $V\in\g{g}_\alpha$. The Levi-Civita connection of $AN$ is
given by
\[
{\nabla}_{aB+U+xZ}(bB+V+yZ)
=\left(\frac{1}{2}\langle U,V\rangle+xy\right)\!B-\frac{1}{2}\left(bU+yJU+xJV\right)
+\left(\frac{1}{2}\langle JU,V\rangle-bx\right)\!Z,
\]
where $a$, $b$, $x$, $y\in\R$, $U$, $V\in\g{g}_\alpha$, and all
vector fields are considered to be left-invariant.

\section{The examples}\label{secExamples}

In this section we present examples of isoparametric hypersurfaces in
$\C H^n$, $n\geq 2$, that are in general not homogeneous. These
hypersurfaces will be tubes around certain homogeneous submanifolds,
so we proceed first with the construction of the latter.

Let $\g{w}$ be a subspace of $\g{g}_\alpha$ and define
$\g{w}^\perp=\g{g}_\alpha\ominus\g{w}$, the orthogonal complement of
$\g{w}$ in $\g{g}_\alpha$. We also define $k=\dim\g{w}^\perp$. Then,
$\g{s}_{\g{w}}=\g{a}\oplus\g{w}\oplus\g{g}_{2\alpha}$ is a solvable
Lie subalgebra of $\g{a}\oplus\g{n}$, as one can easily check from
the bracket relations above. Let $S_{\g{w}}$ be the corresponding
connected subgroup of $AN$ whose Lie algebra is $\g{s}_{\g{w}}$. We
define the submanifold $W_{\g{w}}$ as the orbit $S_{\g{w}}\cdot o$ of
the group $S_{\g{w}}$ through the point $o$. Hence, $W_{\g{w}}$ is a
homogeneous submanifold of $\C H^n$ of codimension $k$.


For each $\xi\in\g{w}^\perp$, we write $J\xi= P\xi+F\xi$, where
$P\xi$ is the orthogonal projection of $J\xi$ onto $\g{w}$, and
$F\xi$ is the orthogonal projection of $J\xi$ onto $\g{w}^\perp$. We
define the K\"{a}hler angle of $\xi\in\g{w}^\perp$ with respect to
$\g{w}^\perp$ as the angle $\varphi_\xi\in\left[0,\pi/2\right]$
between $J\xi$ and $\g{w}^\perp$; hence $\varphi_\xi$ satisfies
$\langle F\xi,F\xi\rangle =(\cos^2\varphi_\xi)\langle
\xi,\xi\rangle$. It readily follows from $J^2=-I$ that $\langle
P\xi,P\xi\rangle =(\sin^2\varphi_\xi)\langle \xi,\xi\rangle$. Hence,
if $\xi$ has unit length, $\varphi_\xi$ is determined by the fact
that $\cos\varphi_\xi$ is the length of the orthogonal projection of
$J\xi$ onto $\g{w}^\perp$. For $\xi\in\g{w}^\perp$ it is convenient
to define
\[
\bar{P}\xi=\frac{1}{\sin\varphi_\xi}P\xi,\text{ (if $P\xi\neq 0$) },
\quad\text{ and }\quad
\bar{F}\xi=\frac{1}{\cos\varphi_\xi}F\xi,\text{ (if $F\xi\neq 0$)}.
\]
Thus, if $\xi$ is of unit length, so are $\bar{P}\xi$ and
$\bar{F}\xi$ if they exist.

Let $\g{c}$ be the maximal complex subspace of $\g{s}_\g{w}$.
Clearly, 
$\g{c}=\g{a}\oplus(\g{g}_\alpha\ominus\C\g{w}^\perp)\oplus\g{g}_{2\alpha}$,
and since $\C\g{w}^\perp=\g{w}^\perp+ J\g{w}^\perp=\g{w}^\perp+
P\g{w}^\perp+ F\g{w}^\perp=\g{w}^\perp\oplus P\g{w}^\perp$, we have
the vector space direct sum decompositions $\g{s}_{\g{w}}=\g{c}\oplus
P\g{w}^\perp$ and $\g{a}\oplus\g{n}=\g{c}\oplus
P\g{w}^\perp\oplus\g{w}^\perp$. Denoting by $\g{C}$, $P\g{W}^\perp$,
and $\g{W}^\perp$ the corresponding left-invariant distributions on
$AN$, we get the tangent bundle $TW_{\g{w}}=\g{C}\oplus P\g{W}^\perp$
and the normal bundle $\nu W_{\g{w}}=\g{W}^\perp$ along $W_{\g{w}}$.

Recall that the shape operator $\mathcal{S}_\xi$ of $W_{\g{w}}$ with
respect to a unit normal $\xi\in\nu W_{\g{w}}$ is defined by
$\mathcal{S}_\xi X=-(\nabla_X\xi)^\top$, for any $X\in TW_{\g{w}}$,
and where $(\cdot)^\top$ denotes orthogonal projection onto the
tangent space. The expression for the Levi-Civita connection of $AN$
allows us to calculate the shape operator  of $W_\g{w}$ for
left-invariant vector fields:
\begin{align*}
\mathcal{S}_\xi U & =0,\quad\text{ if $U\in\g{c}\ominus\g{g}_{2\alpha}$},\qquad
\mathcal{S}_\xi Z =\frac{1}{2}P\xi=\frac{1}{2}(\sin\varphi_\xi)\bar{P}\xi,\\
\mathcal{S}_\xi \bar{P}\eta & =-\frac{1}{2}\langle J\bar{P}\eta,\xi\rangle Z
    =\frac{1}{2}\langle \bar{P}\eta, P\xi\rangle Z= 0,
    \quad\text{ if $\bar{P}\eta\in P\g{w}^\perp\ominus\R P\xi$},\\
\mathcal{S}_\xi \bar{P}\xi & = -\frac{1}{2}\langle J\bar{P}\xi,\xi\rangle Z
    =\frac{1}{2}\langle \bar{P}\xi, P\xi\rangle Z= \frac{1}{2}(\sin\varphi_\xi) Z.
\end{align*}
Note that, if $\varphi_\xi=0$ (that is, if $P\xi=0$), then
$\mathcal{S}_\xi=0$. If $\varphi_\xi>0$ then $W_\g{w}$ has three
principal curvatures with respect to the unit normal vector $\xi$
\[
\frac{1}{2}\sin\varphi_\xi,\quad -\frac{1}{2}\sin\varphi_\xi,\quad\text{and}\quad 0,
\]
whose principal spaces are $\R(Z+\bar{P}\xi)$, $\R(Z-\bar{P}\xi)$,
and $\g{a}\oplus(\g{w}\ominus \R \bar{P}\xi)$, respectively. In any
case, the submanifold $W_\g{w}$ is minimal.\medskip

Now we construct isoparametric hypersurfaces in $\C H^n$. Denote by
$M^r$ the tube of radius $r$ around the submanifold $W_\g{w}$. We
claim that, for every $r>0$, $M^r$ is an isoparametric real
hypersurface which has, in general, nonconstant principal curvatures.
We may assume $k>1$, as otherwise we just obtain parallel
hypersurfaces to the hypersurface $W^{2n-1}$ studied by Berndt
in~\cite{B98}.

We will use Jacobi field theory to calculate the mean curvature of
$M^r$. Given a geodesic $\gamma$ in $\C H^n$, a field $\zeta$ along
$\gamma$ satisfies the Jacobi equation in $\C H^n$ if
$4\zeta''-\zeta-3\langle\zeta,J\dot{\gamma}\rangle J\dot{\gamma}=0$.
Let $p\in W_\g{w}$, and $\xi\in\nu_p W_\g{w}$ be a unit normal
vector. Denote by $\gamma_\xi$ the geodesic such that
$\gamma_\xi(0)=p$ and $\dot{\gamma}_\xi(0)=\xi$. We denote by
$\mathcal{B}_v(t)$ the parallel translation of $v\in T_p\C H^n$ along
$\gamma_\xi$ from $p=\gamma_\xi(0)$ to $\gamma_\xi(t)$. If $X\in T_p
W_\g{w}$ is such that $\mathcal{S}_\xi X=\lambda X$, the solution
$\zeta_X$ to the Jacobi equation with initial conditions
$\zeta_X(0)=X$ and $\zeta'_X(0)=-\mathcal{S}_\xi X$ is given by
\[
\zeta_X(t)=f_\lambda(t) \mathcal{B}_X(t)
    +\langle X, J\xi\rangle g_\lambda(t) J\dot{\gamma}_\xi(t),
\]
where
\[
f_\lambda(t)=\cosh\frac{t}{2}-2\lambda\sinh\frac{t}{2},\qquad
g_\lambda(t)=\left(\cosh\frac{t}{2}-1\right)
\left(1+2\cosh\frac{t}{2}-2\lambda\sinh\frac{t}{2}\right).
\]
On the other hand, for every $\eta\in \nu_p W_\g{w}\ominus \R \xi$,
the solution $\zeta_\eta$ to the Jacobi equation with initial
conditions $\zeta_\eta(0)=0$ and $\zeta'_\eta(0)=\eta$ is given by
\[
\zeta_\eta(t)=p(t) \mathcal{B}_\eta(t)+\langle \eta, J\xi\rangle q(t) J\dot{\gamma}_\xi(t),
\]
where
\[
p(t)=2\sinh\frac{t}{2},\qquad
q(t)=2\sinh\frac{t}{2}\left(\cosh\frac{t}{2}-1\right).
\]

Hence, the above expression for the shape operator $\mathcal{S}_\xi$
of $W_\g{w}$ with respect to a unit normal $\xi$, allows us to
compute (we give the explicit calculations for
$\varphi_\xi\in(0,\pi/2)$; notice that some adaptations are needed in
case ${P}\xi=0$ or ${F}\xi=0$):
\begin{align*}
\zeta_X(t)={}&\cosh\frac{t}{2}\,B_X(t),
    \quad\text{ if $X\in T W_\g{w}\ominus(\R Z\oplus\R\bar{P}\xi)$},\\
\zeta_Z(t)={}&\cosh\frac{t}{2}\,B_Z(t)
    -\sin\varphi_\xi\left(\cos^2\varphi_\xi
    +\sin^2\varphi_\xi\cosh\frac{t}{2}\right)\sinh\frac{t}{2}\,
    \mathcal{B}_{\bar{P}\xi}(t)\\
    &{}-\cos\varphi_\xi\sin^2\varphi_\xi
    \left(\cosh\frac{t}{2}-1\right)
    \sinh\frac{t}{2}\,\mathcal{B}_{\bar{F}\xi}(t),\\
\zeta_{\bar{P}\xi}(t)={}&{}-\sin\varphi_\xi\sinh\frac{t}{2}\,\mathcal{B}_Z(t)
    +\left(\cos^2\varphi_\xi\cosh\frac{t}{2}+\sin^2\varphi_\xi
    \cosh t\right)\mathcal{B}_{\bar{P}\xi}(t)\\
    &{}-\sin\varphi_\xi\cos\varphi_\xi
    \left(\cosh\frac{t}{2}-\cosh t\right)\mathcal{B}_{\bar{F}\xi}(t),\\
\zeta_{\bar{F}\xi}(t)={}&2\sin\varphi_\xi\cos\varphi_\xi
    \left(\cosh\frac{t}{2}-1\right)\sinh\frac{t}{2}\,\mathcal{B}_{\bar{P}\xi}(t)\\
    &{}+2\left(1+\cos^2\varphi_\xi\left(\cosh\frac{t}{2}-1\right)\right)
    \sinh\frac{t}{2}\,\mathcal{B}_{\bar{F}\xi}(t),\\
\zeta_X(t)={}&2\sinh\frac{t}{2}\,\mathcal{B}_X(t),\quad
    \text{ if  $X\in\nu W_\g{w}\ominus(\R\xi\oplus\R\bar{F}\xi)$}.
\end{align*}

We define the endomorphism $D(r)$ of $T_{\gamma_{\xi}(r)}M^r$ by
$D(r)\mathcal{B}_X(r)=\zeta_X(r)$ for each $X\in T_p\C H^n\ominus \R
\xi$. If $\varphi_\xi=0$ (that is, ${P}\xi=0$) then $D(r)$ can be
represented by a diagonal matrix with blocks $\cosh(t/2)I_{2n-k}$,
$2\sinh(t/2)I_{k-2}$, and $2\sinh(t/2)\cosh(t/2)I_1$, the last one
corresponding to the vector $\bar{F}\xi=J\xi$. If $\varphi_\xi=\pi/2$
(that is, ${F}\xi=0$) then $D(r)$ has two diagonal blocks
$\cosh(t/2)I_{2n-k-2}$ and $2\sinh(t/2)I_{k-1}$, and a $2\times 2$
block corresponding to the vectors $Z$, and $\bar{P}\xi=J\xi$. If
$\varphi_\xi\in(0,\pi/2)$ then $D(r)$ has two diagonal blocks
$\cosh(t/2)I_{2n-k-2}$ and $2\sinh(t/2)I_{k-2}$, and another $3\times
3$ block corresponding to the vectors $Z$, $\bar{P}\xi$, and
$\bar{F}\xi$.

It is well known that if $D(r)$ is nonsingular for each $\xi\in\nu
W_{\g{w}}$, then $M^r$ is a hypersurface of $\C H^n$. In our case,
regardless of the value of $\varphi_\xi$, we have
\[
\det(D(r))=2^{k-1}\left(\cosh\frac{r}{2}\right)^{2n-k+1}\left(\sinh\frac{r}{2}\right)^{k-1},
\]
and hence, $M^r$ is a hypersurface for every $r>0$. In this
situation, Jacobi field theory shows that the shape operator of $M^r$
at $\gamma_\xi(r)$ with respect to $-\dot{\gamma}_\xi(r)$ is given by
$\mathcal{S}^r=D'(r)D(r)^{-1}$. Therefore, the mean curvature
$\mathcal{H}^r$ of $M^r$ is
\[
\mathcal{H}^r(\gamma_\xi(r))=\tr\mathcal{S}^r(\gamma_\xi(r))
=\frac{\frac{d}{dr}\det(D(r))}{\det(D(r))}
=\frac{1}{2\sinh\frac{r}{2}\cosh\frac{r}{2}} \left(k-1+2n\sinh^2\frac{r}{2}\right).
\]
Notice again that this value does not depend on the unit vector
$\xi\in\nu W_{\g{w}}$. Therefore, for every $r>0$, the tube $M^r$
around $W_\g{w}$ is a hypersurface with constant mean curvature, and
hence, tubes around the submanifold $W_\g{w}$ constitute an
isoparametric family of hypersurfaces in $\C H^n$, that is, every
tube $M^r$ is an isoparametric hypersurface.

Moreover, it was proved in~\cite{BB01} that the tubes around
$W_{\g{w}}$ are homogeneous precisely when $\g{w}^\perp$ has constant
K\"{a}hler angle, that is, when $\varphi_\xi$ is independent of the
vector $\xi\in\g{w}^\perp$. Indeed, the Berndt-Br\"{u}ck submanifolds
$W_\varphi^{2n-k}$~\cite{BD09} are precisely those $W_{\g{w}}$ for
which $\g{w}^\perp$ has \emph{constant} K\"{a}hler angle $\varphi$ and
$k=\dim\g{w}^\perp$. We summarize all this in the following

\begin{theorem}\label{thMain}
Let $\g{g}=\g{k}\oplus\g{p}$ be the Cartan decomposition of the Lie
algebra of the isometry group $G=SU(1,n)$ of $\C H^n$ with respect to
a point $o\in\C H^n$. Assume $\g{a}\subset\g{p}$ is a maximal abelian
subspace and let
$\g{g}=\g{g}_{-2\alpha}\oplus\g{g}_{-\alpha}\oplus\g{g}_{0}
\oplus\g{g}_{\alpha}\oplus\g{g}_{2\alpha}$ be the root space
decomposition with respect to $\g{a}$. Let $W_{\g{w}}$ be the orbit
through $o$ of the connected subgroup $S_{\g{w}}$ of $G$ whose Lie
algebra is $\g{s}_{\g{w}}=\g{a}\oplus\g{w}\oplus\g{g}_{2\alpha}$,
where $\g{w}$ is any proper subspace of $\g{g}_\alpha$.

Then, the tubes around the submanifold $W_{\g{w}}$ are isoparametric
hypersurfaces of $\C H^n$, and are homogeneous if and only if
$\g{w}^\perp=\g{g}_\alpha\ominus\g{w}$ has constant K\"{a}hler angle.
\end{theorem}

It is feasible to calculate the shape operator $\mathcal{S}^r$ of the
tube $M^r$ at $\gamma_\xi(r)$ from the formula
$\mathcal{S}^r=D'(r)D(r)^{-1}$, although calculations are very long.
We do not give $\mathcal{S}^r$ but mention that its characteristic
polynomial is
\[
p_{r,\xi}(x)=(\lambda-x)^{2n-k-2}\left(\frac{1}{4\lambda}-x\right)^{k-2}q_{r,\xi}(x),
\]
where $\lambda=\frac{1}{2}\tanh\frac{r}{2}$ and
\[
q_{r,\xi}(x)=-x^3+\left(3\lambda+\frac{1}{4\lambda}\right)x^2
-\frac{1}{2}\left(6\lambda^2+1\right)x
+\frac{16\lambda^4+16\lambda^2-1+(4\lambda^2-1)^2\cos 2\varphi_\xi}{32\lambda}.
\]
(This is the same as~\cite[p. 146]{BD09}.)

It is important to remark that, \emph{pointwise}, $M^r$ has the same
principal curvatures as the tubes around the Berndt-Br\"uck
submanifolds $W^{2n-k}_\varphi$, $\varphi\in[0,\pi/2]$,
$k\in\{1,\dots,n-1\}$ (see~\cite{BD09}); notice that for $\varphi=0$
these are tubes around a totally geodesic $\C H^k$,
$k\in\{1,\dots,n-1\}$ in $\C H^n$. In other words, at each point, the
tubes around $W_{\g{w}}$ have the same principal curvatures (with the
same multiplicities) as the homogeneous hypersurfaces that arise as
tubes around the $W_\varphi^{2n-k}$. However, in general, the
principal curvatures, and even the number of principal curvatures,
vary from point to point in $M^r$. Again, \emph{the principal
curvatures of $M^r$ are constant if and only if $\g{w}^\perp$ has
constant K\"{a}hler angle}, that is, if $\varphi_\xi$ does not depend on
$\xi$; this corresponds precisely to the homogeneous examples
constructed by Berndt and Br\"{u}ck~\cite{BB01}.

If $n=2$, then either $\g{w}^\perp$ is $1$-dimensional, in which case
$W_{\g{w}}$ is the Lohnherr hypersurface $W^{2n-1}$, whose
equidistant hypersurfaces are homogeneous~\cite{B98}, or
$\g{w}^\perp=\g{g}_\alpha$, which gives a totally geodesic $\C H^1$
and thus the tubes around it are also homogeneous. In any case, for
$n=2$ we do not get inhomogeneous examples. However, for $n\geq 3$
this construction yields inhomogeneous hypersurfaces in $\C H^n$, for
appropriate choices of $\g{w}$.

%
%
%
%


\end{document}